\documentclass[12pt,bezier]{article}
\usepackage{amsmath, amssymb, amsfonts, amsthm,tikz, url,xcolor}

\newtheorem{thm}{Theorem}

\newtheorem{lem}[thm]{Lemma}
\newtheorem{rem}[thm]{Remark}
\newtheorem{conj}[thm]{Conjecture}

\newcommand{\0}{{\bf 0}}
\newcommand{\1}{{\bf 1}}
\newcommand{\bfa}{{\bf a}}
\newcommand{\bfb}{{\bf b}}
\newcommand{\bfc}{{\bf c}}

\newcommand{\bfv}{{\bf v}}
\newcommand{\x}{{\bf x}}
\newcommand{\y}{{\bf y}}

\newcommand{\corank}{{\rm corank}}
\newcommand{\rank}{{\rm rank}}
\newcommand{\col}{{\rm Col}}
\newcommand{\row}{{\rm Row}}
\newcommand{\wt}{{\rm wt}}
\newcommand{\B}{{\cal B}}
\newcommand{\BB}{{\mathbb B}}

\newcommand{\om}{\Omega}
\newcommand{\ov}{\overline}
\newcommand{\e}{\epsilon}
\newcommand\blfootnote[1]{
	\begingroup
	\renewcommand\thefootnote{}\footnote{#1}
	\addtocounter{footnote}{-1}
	\endgroup}

\title{A Discrete Variation of Littlewood--Offord Problem}

\author{Hossein Esmailian$^{\,\rm a}$\qquad Ebrahim Ghorbani$^{\,\rm a,b,}$\thanks{Corresponding author}
  \\[.4cm]
	{\sl $^{\rm a}$Department of Mathematics, K. N. Toosi University of Technology,}\\
	{\sl P. O. Box 16765-3381, Tehran, Iran}\\
	{\sl $^{\rm b}$Department of Mathematics, University of Hamburg, }\\
	{\sl Bundesstra\ss e 55 (Geomatikum), 20146 Hamburg, Germany }}

\begin{document}
\maketitle
\blfootnote{{\em E-mail Addresses}: {\tt h.esmailian@ipm.ir} (H. Esmailian), {\tt e\_ghorbani@ipm.ir} (E. Ghorbani) }

\begin{abstract}
Littlewood--Offord Problem concerns the number of subsums of a set of vectors that fall in a given convex set. We present a discrete variation of this problem where we estimate the number of subsums that are $(0,1)$-vectors.
 We then utilize this to find the maximum order of graphs with given rank or corank. The rank of a graph $G$  is the rank of its adjacency matrix $A(G)$ and the corank of $G$ is  the rank of $A(G)+I$.

\vspace{4mm}

\noindent {\bf AMS Mathematics Subject Classification\,(2010):}   05C50, 05C75, 15A03.
\end{abstract}

\section{Introduction}

\subsection{Littlewood--Offord Problem and its variants }
Littlewood and Offord \cite{lo} dealt with the following problem in studying the number
of real zeros of random polynomials:  given $\ell$ complex numbers
of modulus at least $1$, from all $2^\ell$ subsums, at
most how many can differ from each other by less than $1$?
They obtained  the bound ${\cal O}\left(\frac{\log \ell}{\sqrt \ell}2^\ell\right)$, which was good enough for
their purpose.
Erd\H os \cite{erd} noticed that for real numbers, Sperner's theorem (stating that any family of subsets of an $\ell$-set no two of which being comparable by inclusion has size at most ${\ell\choose\lfloor\frac \ell2\rfloor}$)
implies a best possible bound.
Suppose $x_1,\ldots,x_\ell$ are real
numbers of modulus at least $1$. For $S\subset\{1,\ldots,\ell\}$, set $x_S =\sum_{i\in S}x_i$.
Then $|x_S-x_{S'}|<1$ implies that $S$ and $S'$ are not comparable by inclusion. So Sperner's theorem implies the following:
\begin{thm}[Erd\H os \cite{erd}] Let $x_1,\ldots,x_\ell$ be real numbers with $|x_i|\ge1$ for all $i$. Let
$\Lambda$ be an open interval of length $1$. Then the total
number of $\ell$-tuples $(\e_1, \ldots,\e_\ell) \in \{0, 1\}^\ell$ with $\e_1x_1 + \cdots+\e_\ell x_\ell\in \Lambda$ is at most ${\ell\choose\lfloor\frac \ell2\rfloor}$.
\end{thm}
This bound is
clearly best possible: if $x_1 = \cdots =x_\ell= 1$, then ${\ell\choose\lfloor\frac\ell2\rfloor}$ of the subsums
are equal to $\lfloor\frac\ell2\rfloor$.
Kleitman \cite{kl65} and Katona \cite{kat} proved that the bound ${\ell\choose\lfloor\frac\ell2\rfloor}$
holds for sums of complex numbers as well.
Later, Kleitman (settling a conjecture of Erd\H os \cite{erd})
proved that instead of complex numbers, vectors in a Hilbert space can be taken.
\begin{thm}[Kleitman \cite{kl70}] Let $\x_1,\ldots,\x_\ell$ be vectors in a Hilbert space, each with
length at least $1$.
Let $\Lambda$ be an open ball of diameter $1$. Then the total
number of $\ell$-tuples $(\e_1, \ldots,\e_\ell) \in \{0, 1\}^\ell$ with $\e_1\x_1 + \cdots+\e_\ell\x_\ell\in \Lambda$ is at most ${\ell\choose\lfloor\frac \ell2\rfloor}$.
\end{thm}

These results attracted the attention of many researchers and
numerous variants of the Littlewood--Offord problem have been proposed and investigated.
 Tao and Vu \cite{TaoVu} initiated a line of work known as inverse Littlewood--Offord theorems.  This theory and its variants played a key role in
estimating the singularity probability of random matrices (see, for instance, \cite{fgg, RudVer,TaoVu,tikh}).
The Littlewood--Offord type theorems has also arisen in
other contexts. In \cite{grig,vw} a modular version of the Littlewood--Offord problem is considered with application to database security.

In the present paper, we address the following question:

\vspace{.25cm}\noindent{\bf Discrete Variation of Littlewood--Offord Problem.} Given $\x_1,\ldots,\x_\ell\in\mathbb{R}^k$, consider all $2^\ell$ subsums $\x_S = \sum_{i\in S}\x_i$ for $S\subseteq\{1,\ldots,\ell\}$. How many of these are $(0,1)$-vectors? In other words, among the $2^\ell$ linear combinations of the columns of the matrix $\left[\begin{array}{c|c|c|c}\x_1&\x_2&\cdots&\x_\ell\end{array}\right]$ with $0,1$ coefficients, how many result in a $(0,1)$-vector?

\vspace{.25cm}
We observe that (see Remark~\ref{remark} below) it is enough to consider reduced matrices, that is, matrices with all distinct rows, each having at least two non-zero components. Throughout, we denote the set of all $(0,1)$-vectors of length $\ell$ by $\{0,1\}^\ell$.

\begin{thm}\label{thm:MainIntroduction} Let $\x_1,\ldots,\x_\ell\in\mathbb{R}^k$ such that the matrix $\left[\begin{array}{c|c|c|c}\x_1&\x_2&\cdots&\x_\ell\end{array}\right]$ is reduced. Let $\Lambda=\{0,1\}^k$. Then the total
	number of $\ell$-tuples $(\e_1, \ldots,\e_\ell) \in \{0, 1\}^\ell$ with $\e_1\x_1 + \cdots+\e_\ell\x_\ell\in\Lambda$ is at most 	 $\frac{2^k+1}{2^{k+1}}\cdot 2^\ell$ if $k \leq \ell-1$, and 	$2^{\ell-1}$ if $k\geq \ell$.
\end{thm}
The proof of Theorem~\ref{thm:MainIntroduction} as well as the characterization of the equality cases for $1\leq k \leq \ell-1$ will be given in Section~\ref{sec:litt}. Our subsequent results in Section~\ref{sec:appl} delve into applications of Theorem~\ref{thm:MainIntroduction} within the subject of rank-order problems in algebraic graph theory, a domain that is described in the next subsection.

\subsection{Rank-order problems for graphs}
Let $G$ be a simple graph with vertex set $\{v_1, \ldots, v_n\}$.  The {\em adjacency matrix} of $G$ is an $n \times n $ matrix $A(G)$ whose $(i, j)$-entry is $1$
if $v_i$ is adjacent to $v_j$ and $0$ otherwise.
The {\em order} of $G$ is the number of vertices of $G$. We denote the set of neighbors of a vertex $v$ of $G$ by $N(v)$.
By {\em eigenvalues} and {\em rank} of $G$, we mean the eigenvalues and the rank of $A(G)$ over the reals.  We denote the latter by $\rank(G)$.

Let $\mu$ be a graph eigenvalue.
An extremal problem in algebraic graph theory asks for finding the maximum order $n$ of a graph $G$ where $\rank(A(G)-\mu I)$ is a given integer $r$.
Rowlinson \cite{r1} showed that if $\mu \notin\{0,-1\}$, then $n<r+2^r$. This was improved in \cite{r2} to $n\le\frac{1}{2}r(r+5)-2$.
Bell and Rowlinson \cite{br} finally  proved that if $\mu \notin\{0, -1\}$, then either (i) $n\le \frac{1}{2}r(r+1)$ or (ii) $\mu=1$ and $G=K_2$ or $2K_2$.

As the above result suggests, $\mu =0, -1$ are somewhat exceptional. We first discuss the case of $\mu=0$. In general, the order of graphs $G$ with a fixed $r=\rank(G)$ can be unbounded.
In fact, the order of $G$ can be increased without changing its rank by adding a new vertex $v$ {\em twin} with a vertex $u$
(i.e. with $N(u)=N(v)$) to $G$ or adding isolated vertices.
For this reason, only {\em reduced graphs}, that is,
graphs with no isolated vertices and no twins are taking into account.
For the reduced graphs with rank $r$, it is easily seen that the order is bounded above by $2^r-1$.
This bound is far from being sharp. Kotlov and Lov\'asz \cite{kl} solved the problem asymptotically.
They proved that any reduced graph of rank $r$ has order ${\cal O}(2^{r/2})$ and, for every $r\ge2$, they constructed a reduced graph of rank $r$ and order
$$	m(r)=\left\{
	\begin{array}{ll}
		2^{\frac{r}{2}+1}-2& r\ \text{even,}\\
		5\cdot 2^{\frac{r-3}{2}}-2& r\ \text{odd.}
	\end{array}
	\right.$$
 This is  conjectured to be the precise value of the maximum order:
\begin{conj}[Akbari, Cameron and Khosrovshahi \cite{ack}]\label{conj:ack}
	The maximum order of a reduced graph with rank $r\ge2$ is equal to $m(r)$.
\end{conj}
Haemers and Peeters \cite{rank10} proved Conjecture~\ref{conj:ack} for graphs containing an induced matching of size $r/2$ for even $r$ or an induced
subgraph consisting of a matching of size $(r-3)/2$ and a cycle of length $3$ for odd $r$. Ghorbani, Mohammadian and Tayfeh-Rezaie   \cite{gmt3} proved
that if Conjecture~\ref{conj:ack} is valid for all reduced graphs of rank at most $46$, then it is true in general.
Further, they showed that the order of every reduced graph of rank $r$ is at most $8m(r)+14$. This problem has been also investigated within specific
families of graphs. In \cite{gmt}, it is proved that the maximum order of every reduced tree and  bipartite graph of rank $r$ is
$3r/2-1$ and $2^{r/2}+r/2-1$, respectively. This value is shown to be
$3\cdot 2^{\lfloor {r/2}\rfloor-2}+\lfloor r/2\rfloor$
for non-bipartite triangle-free graphs in \cite{gmt2}.

 For the other exceptional eigenvalue, namely $\mu=-1$, one should consider the rank of $A(G)+I$ which we call it the {\em corank} of
$G$  denoted by $\corank(G)$. Similar to the case of rank, the order of graphs with a fixed corank can be unbounded.
In fact, in any graph $G$, adding a new vertex $v$ {\em cotwin} with a vertex $u$ (i.e. with $N(u)\cup\{u\}=N(v)\cup\{v\}$) to $G$,
increases the order of $G$ without changing its corank. Therefore, one should consider {\em coreduced} graphs, i.e. graphs with no cotwins.
Similar to the case of rank, in \cite{esg}, we showed that the order of coreduced graphs with corank $r$ is ${\cal O}(2^{r/2})$.
It was also shown that the order of any tree and bipartite
graph of corank $r$ is at most $2r-3$ and  $2r-2$, respectively, and  the order of any coreduced cotree (i.e. the complement of a tree) of $\corank$ $r$ is at most $\left\lfloor3r/2-2\right\rfloor$.

As applications for our discrete variation of Littlewood--Offord Problem, we (i) determine the maximum order of a coreduced graph with a bipartite complement of given corank, and (ii) give a new proof for the result of \cite{gmt} on the maximum order of a reduced bipartite graph of given rank. In both cases, we characterize the graphs achieving the maximum order. These results will be presented in Section~\ref{sec:appl}.

\section{Discrete Variation of Littlewood--Offord Problem}\label{sec:litt}
Our objective in this section is to prove Theorem~\ref{thm:MainIntroduction}. Some notation is in order.
In the remainder of the paper all vectors are treated as ``row vectors.''
Let $\bfv$ be a real vector. The {\em weight} of $\bfv$, denoted by $\wt(\bfv)$, is the number of non-zero components of $\bfv$.
Let $A$ be a $k\times \ell$ matrix.
We set$$\om(A) := \big\{\bfb \in \{0,1\}^\ell : \bfb A^\top\in \{0,1\}^k\big\}.$$
In other words, $\om(A)$ is the set of $(0,1)$-vectors $\bfb$ of length $\ell$ such that the linear combination of the columns of $A$ with the coefficients from $\bfb$ gives a $(0,1)$-vector.
As a discrete variation of Littlewood--Offord Problem, in this section we deal with estimating the size of $\om(A)$.
We call a real matrix {\em reduced} if all its rows are distinct and have weight at least $2$. Our main result is that if $A$ is reduced, then $\om(A)$ has size at most  $2^{\ell-1}$ for $k\geq \ell$, and
$\frac{2^k+1}{2^{k+1}}\cdot 2^\ell$ for $k \leq \ell-1$.

\begin{rem}\label{remark}\rm
Here we justify the restriction to the reduced matrices.
 If $\bfv$ is vector of length $\ell$, then
$\om(\bfv)$ is the set of all $\bfb \in \{0,1\}^\ell$ such that the inner product $\bfv\cdot\bfb$ is $0$ or $1$.
Note that if $\bfv_1,\ldots,\bfv_k$ are all the rows of $A$, then
\begin{equation}\label{eq:omA,omvk}
\om(A)=\om(\bfv_1)\cap\cdots\cap\om(\bfv_k).
\end{equation}
So deleting repeated rows does not alter $\om(A)$.
If some $\bfv_i$ has weight one and its non-zero component is not $1$, then
$|\om(\bfv_i)|=2^{\ell-1}$, and thus by \eqref{eq:omA,omvk}, $|\om(A)|\le2^{\ell-1}$, so we are done.
Otherwise,
assume that any weight-one row $\bfv_i$ is a $(0,1)$-vector. In that case,  $\om(\bfv_i)=\{0,1\}^\ell$.
It follows that $\om(A)=\om(A')$ where $A'$ is obtained from $A$ be removing repeated rows as well as any row of weight at most $1$.
\end{rem}

As we shall see, our main problem on bounding
$|\om(A)|$ for real matrices $A$, can be reduced to $(0,\pm1)$-matrices.
So in the next few lemmas, we deal with matrices/vectors with $0,\pm1$ entries.

\begin{lem}\label{lem:Omega(v)}
Let $\bfv$ be a $\pm1$-vector of length $\ell$. If the number of $1$'s in $\bfv$ is $k$,
then, $|\Omega(\bfv)|= {\ell+1 \choose k}\le {\ell+1 \choose \lfloor\frac{\ell+1}{2}\rfloor}$.
\end{lem}
\begin{proof}
With no loss of generality, we may assume that $\bfv=(1,\ldots,1,-1,\ldots,-1),$ where the number of $1$'s is $k$.
Let $\bfb=(b_1, \ldots, b_\ell)\in\Omega(\bfv)$ and $\bfb'=(1-b_1, \ldots, 1-b_k, b_{k+1}, \ldots, b_\ell)$. Assume that $\wt((b_1, \ldots, b_k))=s$ and  $\wt((b_{k+1}, \ldots, b_\ell))=t$. Hence $\wt(\bfb')=k-s+t$. We have $s-t=\bfb\cdot\bfv \in \{0,1\}$ and hence $\wt(\bfb')\in \{k, k-1\}$.
 So the number of different $\bfb'$ (and so the number of different $\bfb\in\Omega(\bfv)$) is equal to ${\ell \choose k-1}+{\ell \choose k}={\ell+1 \choose k}$.
  We know that ${\ell+1 \choose k}\le {\ell+1 \choose \lfloor\frac{\ell+1}{2}\rfloor}$, so the proof is complete.
  \end{proof}

Given a matrix $A$, we denote its submatrix consisting of all the non-zero columns  by $A^*$. If $A^*$ is obtained by removing $j$ zero columns, then it is clear that
\begin{equation}\label{eq:A*}
|\om(A)|=2^j\cdot|\om(A^*)|.
\end{equation}
We say that the matrix $A'$ is {\em equivalent} with $A$ and write $A'\simeq A$, if $A$ can be transformed into $A'$ by row and/or column permutations.
It is observed that
$$|\om(A')|=|\om(A)|.$$
From \eqref{eq:omA,omvk}, it is also clear that if the matrix $B$ is obtained by removing some of the rows of $A$, then
$$|\om(A)|\le|\om(B)|.$$
We denote the all $1$'s and all $0$'s vectors by $\1$ and $\0$, respectively.

\begin{lem}\label{lem:2^k(3)}
Let $A$ be a $k\times (k+2)$ matrix of the form
\begin{equation}\label{eq:k=3}
\left[
\begin{array}{cc|cccc}
\pm1&\pm1&\pm1&\ldots&0&0\\
\vdots&\vdots&\vdots&\ddots&\vdots&\vdots\\
\pm1&\pm1&0&\ldots&\pm1&0\\
\pm1&\pm1&0&\ldots&0&a
\end{array}\right],
\end{equation}
where $ a \in\{0,\pm1\}$.  Then $|\om(A)|\le2^{k+1}+2$ and the equality holds if and only if $A$ is of the form
\begin{equation}\label{eq:A1A2}
A_1 = \left[
\begin{array}{cc|c|c}
1&1&&0\\
\vdots&\vdots&-I_{k-1}&\vdots\\
1&1&&0\\ \hline
1&1&{\bf0}&b
\end{array}\right],
\quad
A_2 = \left[
\begin{array}{cc|c|c}
a_1&-a_1&&0\\
\vdots&\vdots&I_{k-1}&\vdots\\
a_{k-1}&-a_{k-1}&&0\\ \hline
1&-1&{\bf0}&c
\end{array}\right],
\end{equation}
where  $a_i \in\{1,-1\}$,  $b \in\{0,-1\}$ and $c \in\{0,1\}$.
\end{lem}
\begin{proof}
If in some row of $A$ with weight $3$ there are not two $1$'s, then by Lemma~\ref{lem:Omega(v)} and \eqref{eq:A*}, $|\om(A)| \le{4\choose1}\cdot2^{k-1}=2^{k+1}$ and we are done. So assume that in any row of $A$ with weight 3, there are exactly two $1$'s.
First, suppose that in the right block of $A$ there exist two entries with different signs.
 Then $A$ contains a $2\times(k+2)$ submatrix $B$ with
$$B^* =\left[
\begin{array}{crrr}
1&-1&1&0\\
1&1&0&-1
\end{array}\right].$$
We see that
$$\om(B^*) =\{ 0000, 0010, 0110, 0111 , 1000, 1001, 1101, 1111\}.$$
Thus $|\om(A)|\le|\om(B)| = |\om(B^*)|\cdot 2^{k-2}= 2^{k+1}$, and so we are done.
Hence, we assume that in the right block of $A$ all the non-zero entries have the same sign.
It follows that  $A$ is of the form either $A_1$ or $A_2$.
We have
$$\om(A_1)=\begin{cases}
\{\0,\01\}\cup\left(\{01,10\}\times\{0,1\}^k\right) & \hbox{if $b=0$},\\
\{\0,\1\}\cup\left(\{01,10\}\times\{0,1\}^k\right) & \hbox{if $b=-1$}.
\end{cases}$$
For $A_2$, consider the $(0,1)$-vectors
$\bfb=\frac12(1-a_1,\ldots, 1-a_{k-1})$ and $\bfb'=\frac12(1+a_1,\ldots, 1+a_{k-1})$. Then
$$\om(A_2)=
\begin{cases}
\{10\bfb0,10\bfb1\}\cup\left(\{00,11\}\times\{0,1\}^k\right) & \hbox{if $c=0$},\\
\{10\bfb0,01\bfb'1\}\cup\left(\{00,11\}\times\{0,1\}^k\right) & \hbox{if $c=1$}.
\end{cases}$$
Therefore,
$|\om(A_1)| =|\om(A_2)| = 2^{k+1}+2$.
\end{proof}
Similar to Lemma~\ref{lem:2^k(3)}, the following can be obtained.
\begin{lem}\label{lem:2^k(2)}
Let $A$ be $k\times (k+1)$ matrix of the form
\begin{equation}\label{eq:k=2}
\left[\begin{array}{c|ccc}
\pm1&\pm1&\ldots&0\\
\vdots&\vdots&\ddots&\vdots\\
\pm1&0&\ldots&\pm1
\end{array}\right].\end{equation}
Then $|\om(A)| \leq 2^k +1$.
The equality holds if and only if $A$ is one of the following matrices:
\begin{equation}\label{eq:A3A4}
A_3=\left[
\begin{array}{c|c}
1&\\
\vdots&-I_k\\
1&\\
\end{array}\right],
\quad
A_4=\left[
\begin{array}{c|c}
\pm1&\\
\vdots&I_k\\
\pm1&\\
\end{array}\right].
\end{equation}
\end{lem}

We also need the following lemma on $(0,\pm1)$-matrices with two or three rows.
\begin{lem}\label{lem:comput} Let $A$ be a $k\times s$ reduced $(0,\pm1)$-matrix and $t$ be the maximum weight of the rows of $A$.
\begin{itemize}
  \item[\rm(i)] If $k=2$, $t=6,7$ and $s\le14$, then $|\om(A)|\le2^{s-1}$.
\item[\rm(ii)] If $k=2$, $t=4,5$ and $s\le10$, then $|\om(A)|<\frac{5}{8}\cdot2^{s}$.
\item[\rm(iii)] If $k=2$, $t=3$, $s\le6$,  and $A^*$ is not equivalent with
\begin{equation}\label{eq:B0}
B_0=\left[\begin{array}{cccc}
\pm1&\pm1&\pm1&0\\
\pm1&\pm1&0&a
\end{array}\right],
\end{equation}
 where $a\in\{0,\pm1\}$, then $|\om(A)|\le\frac{9}{16}\cdot2^{s}$.
  \item[\rm(iv)]  If $k=3$, $t=4,5$ and $s\le15$, then $|\om(A)|\le2^{s-1}$.
  \item[\rm(v)] If $k=3$, $t=3$, $s\le9$, and $A^*$ is not equivalent with the matrix given in \eqref{eq:k=3}, then $|\om(A)|\le2^{s-1}$.
 \end{itemize}
\end{lem}

We verified Lemma~\ref{lem:comput} by performing an exhaustive computer search.\footnote{The Python code of the program is available at: \url{https://wp.kntu.ac.ir/ghorbani/ComputFiles/PythonCode.txt}}
As it may not be clear from the statement, we discuss here why such a search is feasible.
As an instance, we give an enumeration on the total number of inner products required to verify the part (i) of the lemma with $t=7$.
Let $\bfv$ be the first row of $A$ of weight $7$ and $d$ be the number of $1$'s in $\bfv$.
If $d\ne4$, then by Lemma~\ref{lem:Omega(v)}, $|\om(\bfv^*)|\le{8\choose3}<2^6$ implying that   $|\om(\bfv)|\le|\om(\bfv^*)|\cdot2^{s-7}<2^{s-1}$, and we are done.
So let $d=4$. Then $A$ is equivalent
with a matrix of the form
$$\left[
\begin{array}{cccccccccccccc}
-1&-1&-1&0&0&0&0&0&0&0&1&1&1&1\\
a_1&a_2&a_3&b_1&b_2&b_3&
b_4&b_5&b_6&b_7&c_1&c_2&c_3&c_4
\end{array}\right],$$
where $a_1\le a_2\le a_3$, $b_1\le\cdots\le b_7$ and  $c_1\le\cdots\le c_4$.
Let $\bfa=(a_1,a_2,a_3)$, $\bfb=(b_1,\ldots,b_7)$ and  $\bfc=(c_1,\ldots,c_4)$. We must have $2\le\wt(\bfa)+\wt(\bfb)+\wt(\bfc)\le7$. If $\wt(\bfb)=7$, then $\wt(\bfa)=\wt(\bfc)=0$ and thus
$|\om(A)|=|\om(\bfv^*)|\cdot|\om(\bfb)|\le{8\choose4}^2<2^{13}$, and we are done. So $\wt(\bfb)\le6$.
Suppose that $\wt(\bfa)=i$, $\wt(\bfb)=j$ and $\wt(\bfc)=r$. Given that the components of these vectors are increasing, the numbers of choices for $\bfa,\bfb$, and $\bfc$ are $i+1,j+1$, and $r+1$, respectively. 
We have $0\le i\le3$, $0\le j\le6$, and $0\le k\le4$. Furthermore, since $i+j+k\le7$, we must have $j\le7-i$ and $k\le7-i-j$. Taking into account these conditions on $i,j,k$, it follows that the number of different choices for the second row of $A$ is at most
$$\sum_{i=0}^3(i+1)\sum_{j=0}^{\min(6,7-i)}(j+1)\sum_{r=0}^{\min(4,7-i-j)}(r+1)=1267.$$
Now, for any choice of $A$ we should compute $\x A^\top$ for any $\x\in\{0,1\}^{14}$.
Since $A^*$ has $j+7$ columns, it suffices to compute $\x A^{*\top}$ for any $\x\in\{0,1\}^{j+7}$.
It turns out that the total number of required inner products to verify the assertion is at most $$2\sum_{i=0}^3(i+1)\sum_{j=0}^{\min(6,7-i)}2^{j+7}(j+1)\sum_{r=0}^{\min(4,7-i-j)}(r+1)=3035648,$$
which shows the feasibility of the exhaustive search.

We are now prepared to prove the main result of the paper. For convenience, we repeat Theorem~\ref{thm:MainIntroduction} here, including the equality cases.
\begin{thm}\label{main}
If $A$ is a $k\times\ell$ reduced matrix, then
\[|\om(A)| \leq
\begin{cases}
\frac{2^k+1}{2^{k+1}}\cdot 2^\ell &k \leq \ell-1,\\
2^{\ell-1} &k\geq \ell.
\end{cases}\]
For $1\leq k \leq \ell-1$, the equality holds if and only if $A^*$ is equivalent with one of the matrices
$A_1,A_2,A_3,A_4$ given in \eqref{eq:A1A2} and \eqref{eq:A3A4}.
\end{thm}

\begin{proof}
 We first show that if $A$ has an entry other than $0,\pm1$, then we are done.
To see this, with no loss of generality, assume that  $\bfv =(v_1, v_2, \ldots, v_\ell)$, with $v_1\notin \{0,\pm1\}$,  is some row of $A$.
   Let $\bfa=(1, a_2, \ldots, a_\ell)\in\{0,1\}^\ell$ and $\bfa'=(0, a_2, \ldots, a_\ell)$.
   We claim that at most one of $\bfa$ and $\bfa'$ belong to 
   $\om(\bfv)$, since otherwise
$$|v_1|=|\bfa\cdot \bfv-\bfa'\cdot \bfv|\in\{0,1\},$$
   which is a contradiction. Thus, at most one of $\bfa$ or $\bfa'$ belong to $\om(\bfv)$. This implies that $|\om(A)|\le|\om(\bfv)|\le2^{\ell-1}$.
  So we may assume that all the entries of $A$ are $0,\pm1$.

Assume that the row $\bfv$ with $\wt(\bfv)=t$ has the largest weight among the rows of $A$.
   By Lemma~\ref{lem:Omega(v)} and \eqref{eq:A*},  we have $|\om(\bfv)|\le{t+1 \choose \lfloor\frac{t+1}{2}\rfloor}2^{\ell-t}$. For $t\ge8$, by induction, we have ${t+1 \choose \lfloor\frac{t+1}{2}\rfloor}<2^{t-1}$.
  Hence if $t\ge8$, then $|\om(A)|\le|\om(\bfv)|<2^{\ell-1}$, and we are done. Therefore, we suppose that $t\le7$.
We consider the following four cases.

\vspace{.2cm}\noindent\textbf{Case 1.} $k = 1$

Since $k=1$, and $A$ is a reduced matrix, the weight of each row of $A$ is at least two. Thus, $k\le\ell-1$, which means that we only need to show that $|\om(A)|\le\frac{3}{4}\cdot 2^{\ell}$.

  As $t\ge2$, we have ${t+1 \choose \lfloor\frac{t+1}{2}\rfloor}\le\frac34\cdot2^t$ with equality for $t=2,3$.
Now, from Lemma~\ref{lem:Omega(v)} it follows that
$|\om(A)|\le|\om(\bfv)|\le {t+1 \choose \lfloor\frac{t+1}{2}\rfloor}2^{\ell-t}\le\frac34\cdot 2^{\ell}.$ The equality holds if and only if $t=2,3$ which agrees with the equality cases of the theorem.

\vspace{.2cm}\noindent\textbf{Case 2.} $k = 2$

In this case, we need to show that for $\ell=2$, $|\om(A)|\le2$, and for $\ell\ge3$,  $|\om(A)|\le\frac{5}{8}\cdot 2^{\ell}$. 
(The only possibility for $A$ in the case $\ell=2$ is that $A$ is equivalent to the matrix $B_1$ below.)

First, assume that $t=2$. Then, $A^*$ is equivalent with one of
$$B_1 = \left[
\begin{array}{cc}
\pm1&\pm1 \\
\pm1&\pm1
\end{array}\right],
\quad B_2 = \left[
\begin{array}{cccc}
\pm1&\pm1&0&0 \\
0&0&\pm1&\pm1
\end{array}\right],
\quad B_3 = \left[
\begin{array}{ccc}
\pm1&\pm1&0 \\
\pm1&0&\pm1
\end{array}\right].$$
 It is easy to check that at most  two vectors from $\{0,1\}^2$ can belong to $\om(B_1)$, that is, $|\om(B_1)|\le2$.
   So if $A^*\simeq B_1$, then $|\om(A)|=|\om(A^*)|\cdot2^{\ell-2}\le2^{\ell -1}$, implying the result.
   We have $|\om(B_2)|=|\om((\pm1,\pm1))|^2\le9$. Thus, if $A^*\simeq B_2$, then $|\om(A)|=|\om(B_2)|\cdot2^{\ell-4}=\frac{9}{16} \cdot 2^{\ell}<\frac{5}{8}\cdot 2^{\ell}$, and we are done. Finally, let $A^*\simeq B_3$. By Lemma~\ref{lem:2^k(2)}, $|\om(B_3)|\le 5$. It follows that $|\om(A)|\le\frac{5}{8}\cdot 2^{\ell}$ and the equality holds if and only if $A^*$ is equivalent with $A_3$ or $A_4$ of \eqref{eq:A3A4}.

    If $t=3$, then $A^*$  has $s\le6$ columns because the weight of the second row of $A$ is at most $t$.
    If $A^*$ is not equivalent with $B_0$ of \eqref{eq:B0}, then  Lemma~\ref{lem:comput}\,(iii) implies that $|\om(A^*)|\le\frac{9}{16}\cdot2^{s}$ and thus $|\om(A)|\le\frac{9}{16}\cdot2^{\ell}<\frac58\cdot2^\ell$.
    If $A^*\simeq B_0$,  then $s=4$ and by Lemma~\ref{lem:2^k(3)}, $|\om(A^*)|\le10$. It follows that
 $|\om(A)|=|\om(A^*)|\cdot2^{\ell-4}\leq\frac58\cdot 2^{\ell}$ and the equality holds if and only if $A^*$ is  equivalent with either $A_1$ or $A_2$ of \eqref{eq:A1A2}.

 If $t=4,5$, then $A^*$  has $s\le10$ columns. By Lemma~\ref{lem:comput}\,(ii), $|\om(A^*)|<\frac58\cdot2^s$. It follows that
 $|\om(A)|=|\om(A^*)|\cdot2^{\ell-s}<\frac58\cdot2^\ell$.

 If $t=6,7$, then ${t+1 \choose \lfloor\frac{t+1}{2}\rfloor}=\frac{35}{64}\cdot2^t<\frac58\cdot2^t$.
Then by Lemma~\ref{lem:Omega(v)},
$|\om(A)|\le|\om(\bfv)|\le {t+1 \choose \lfloor\frac{t+1}{2}\rfloor}2^{\ell-t}<\frac58\cdot2^\ell$.

\vspace{.2cm}\noindent\textbf{Case 3.} $k = 3$

In this case, we need to show that for $\ell=2,3$, $|\om(A)|\le2^{\ell-1}$, and for $\ell\ge4$,  $|\om(A)|\le\frac9{16}\cdot 2^{\ell}$.

First, let $t=2$. Comparing the $2\times\ell$ submatrices of $A$ with $B_1,B_2,B_3$ of Case~2, we see that $A$ satisfies in one of the following three cases.
\begin{itemize}
	\item[(i)]	For some $2\times\ell$ submatrix $B$ of $A$, we have $B^*\simeq B_1$. Thus $|\om(A)|\le|\om(B)|\le2^{\ell-1}$.
	\item[(ii)] For all $2\times\ell$ submatrices $B$ of $A$,	
	we have $B^*\simeq B_3$.
	Then
	$A^*$ is equivalent either with the matrix given in \eqref{eq:k=2}, or   with
\begin{equation}\label{eq:k=3,t=2,s=3}
\left[\begin{array}{ccc}
		\pm1&\pm1&0\\
		\pm1&0&\pm1\\
		0&\pm1&\pm1
	\end{array}\right].
\end{equation}
If the former occurs, then by Lemma~\ref{lem:2^k(2)}, $|\om(A)|\leq \frac{2^k+1}{2^{k+1}} \cdot 2^{\ell} = \frac{9}{16}\cdot 2^{\ell}$ and
	the equality holds if and only if $A^*$ is equivalent with $A_3$ or $A_4$ of \eqref{eq:A3A4}.
So assume that $A^*$ is equivalent with
\eqref{eq:k=3,t=2,s=3}.
If some $2\times3$ submatrix $B$ of $A^*$ is equivalent to neither of $A_3,A_4$ of \eqref{eq:A3A4}, then by Lemma~\ref{lem:2^k(2)}, $|\om(A^*)|\le|\om(B)|\le4$. It follows that $|\om(A)|=|\om(A^*)|\cdot2^{\ell-3}\le2^{\ell-1}$, as desired.
Otherwise, $A^*$ is equivalent with either of $$\left[\begin{array}{ccc}
		1&1&0\\
		1&0&1\\
		0&1&1
	\end{array}\right],
	\quad
\left[\begin{array}{ccc}
		1&-1&0\\
		1&0&-1\\
		0&1&1
	\end{array}\right].$$
Then it can be easily checked that $|\om(A^*)|=4$ and thus
 $|\om(A)|\le 4\cdot 2^{\ell-3} = 2^{\ell-1}$, and we are done.	

\item[(iii)] $A$ has two $2\times\ell$ submatrices that are either both equivalent with $B_2$, or one is equivalent with $B_2$ and the other one with $B_3$. It turns out that $A^*$ is equivalent with either of
$$
{\small\left[\begin{array}{cccccc}
\pm1&\pm1&0&0&0&0 \\
0&0&\pm1&\pm1&0&0\\
0&0&0&0&\pm1&\pm1
\end{array}\right],~
\left[
\begin{array}{ccccc}
\pm1&\pm1&0&0&0 \\
0&0&\pm1&\pm1&0\\
0&0&0&\pm1&\pm1
\end{array}\right],~
\left[
\begin{array}{cccc}
\pm1&\pm1&0&0 \\
0&0&\pm1&\pm1\\
0&\pm1&\pm1&0
\end{array}\right].}$$
For the first one, we have $|\om(A^*)|\le|\om((\pm1,\pm1))|^3\le27$, and thus $|\om(A)|\le|\om(A^*)|\cdot2^{\ell-6}<2^{\ell-1}$.
For the second one, $|\om(A^*)|\le|\om((\pm1,\pm1))|\cdot|\om(B_3)|\le15$, and thus $|\om(A)|\le|\om(A^*)|\cdot2^{\ell-5}<2^{\ell-1}$.
For the third one, if we have $|\om(A^*)|\le8$, then it will follow that $|\om(A)|\le2^{\ell-1}$.
Otherwise, $|\om(A^*)|\ge9$. On the other hand, $\om(A^*)\subseteq\om(B_2)$. Since $|\om(B_2)|\le9$, it follows that  $\om(A^*)=\om(B_2)$.
This in turn implies that  $\om(B_2)\subseteq\om(\x)$ where $\x=(0,\pm1,\pm1,0)$. At least one of $0100$ or $1100$ and at least one of $0010$ or $0011$ belong to $\om(B_2)$. This implies that $\x=(0,1,1,0)$. Also $\om(B_2)$ contains a vector of the form $*11*$. Such a vector cannot belong to $\om(\x)$, a contradiction.
\end{itemize}

Next, let $t=3$.
Since the weight of each row of $A$ is at most $t$,  $A^*$  has $s\le9$ columns.
If $A^*$ is not equivalent with the matrix given in \eqref{eq:k=3}, then by Lemma~\ref{lem:comput}\,(v), $|\om(A^*)|\le2^{s-1}$. It follows that $|\om(A)| = |\om(A^*)|\cdot2^{\ell-s}\le2^{\ell-1}$, as desired.
Otherwise, by Lemma~\ref{lem:2^k(3)}, $|\om(A)|\leq \frac{2^k+1}{2^{k+1}} \cdot 2^{\ell} =
\frac{9}{16} \cdot 2^{\ell}$ and
 the equality holds if and only if $A^*$ is equivalent with $A_1$ or $A_2$ of \eqref{eq:A1A2}.

 If $t=4,5$, then $A^*$  has $s\le15$ columns.
 By Lemma~\ref{lem:comput}\,(iv), $|\om(A^*)|\le2^{s-1}$. It follows that
 $|\om(A)|\le|\om(A^*)|\cdot2^{\ell-s}\le2^{\ell-1}$, and we are done.

 If $t=6,7$, in a similar manner as above we are done by Lemma~\ref{lem:comput}\,(i).

\vspace{.2cm}\noindent\textbf{Case 4.} $k\geq 4$

First let $t=2$.
If $A^*$ is equivalent with
  the matrix given in \eqref{eq:k=2}, then by Lemma~\ref{lem:2^k(2)},  $|\om(A)|\leq \frac{2^k+1}{2^{k+1}} \cdot 2^{\ell}$ and
 the equality holds if and only if $A^*$ is equivalent with $A_3$ or $A_4$ of \eqref{eq:A3A4}.
  Otherwise, as shown in Case~3, for some  $3\times\ell$ submatrix $B$ of $A$ we have $|\om(B)|\le2^{\ell-1}$, and so we are done.

If $t=3$, then we are done similarly as for $t=2$.

If $4\le t\le7$, then we are done by Lemma~\ref{lem:comput} as in Case~3.
\end{proof}

\section{Applications}\label{sec:appl}
In this section, we present two applications for our result on the discrete variation of Littlewood--Offord Problem.
We first give a new proof for the result of \cite{gmt} on the maximum order of a reduced bipartite graph with a given rank.
Then we present another application on finding the maximum order of a coreduced {\em cobipartite} graph (i.e. the complement of a bipartite graph) with a given corank.

We need further notation. Let $G$ be a bipartite graph. Then its adjacency matrix can be put in the form:
$$A(G)=
\left[
\begin{array}{c|c}
O&B\\ \hline
B^\top&O
\end{array}
\right].$$
We call $B=B(G)$ a {\em bipartite adjacency matrix} of $G$. When $G$ is connected, this is unique up to permutations of rows and columns.
We denote the $\ell\times2^\ell$ matrix whose columns consist of all $(0,1)$-vectors of length $\ell$ by $\BB_\ell$.
The bipartite graph $G$ with $B(G)=\BB_\ell$ is denoted by $\B_\ell$.  The graph $\B_\ell$ is in fact the {\em incidence graph} of $[\ell]:=\{1,\ldots,\ell\}$ versus ${\cal P}([\ell])$, the power set of $[\ell]$.
We also denote the column space and the row space of a matrix $M$ by $\col(M)$ and $\row(M)$, respectively.

\subsection{Bipartite graphs}

The graph $\B_\ell$ has an isolated vertex. We denote the resulting graph by removing this isolated vertex by $\B'_\ell$.  So $\B'_\ell$ is a reduced bipartite graph of rank $2\ell$ and order $2^\ell + \ell-1$.

As the first application of Theorem~\ref{main}, we give a new proof for the following theorem from \cite{gmt}.
\begin{thm}\label{bipartproof2}
Let $G$ be a reduced bipartite graph  of order $n$ and rank $r$. Then $n \leq 2^{r/2} + r/2 -1$ and the equality holds if and only if $G$ is isomorphic to $\B'_{r/2}$.
\end{thm}

\begin{proof}{ Let $B=B(G)$ be a $p\times q$ matrix with rank $\ell$. We have $r=2\ell$. We can assume that $p\le q$.
 First, suppose that $p=\ell$. Since $G$ is a reduced graph,  $B$ has no two identical columns nor a zero column. Thus $q\le2^\ell-1$ with equality if and only if $B$ is equal to the matrix $\BB_\ell$ whose zero column is removed. It follows that $n=p+q \le 2^\ell + \ell-1$ with equality if and only if  $G$ is isomorphic to $\B'_\ell$.

Now, assume that  $p=\ell+k$ with $k\geq1$.
By performing column-elementary operations, we can find a basis for $\col(B)$ as follows (a permutation of the rows might be also necessary):
$$W =\left[
\begin{array}{c}
I_\ell\\ \hline
C_{k \times \ell}
\end{array}
\right].$$
 Since $G$ is a reduced graph, $W$ has no two identical rows and no zero row. This implies that $C$ is a reduced matrix. Any column of $B$ is a non-zero $(0,1)$-vector, so it is generated by a linear combination of the columns of $W$ if the corresponding vector of coefficients belong to  $\om(W)\setminus\{\0\}$.
 It turns out that $q\le|\om(W)|-1$. It is also clear that $\om(W)=\om(C)$.
 If $k\ge\ell$, by Theorem~\ref{main}, $|\om(C)|\leq 2^{\ell-1}$ and then as $p\le q$, we have
$n=p+q\le 2q \leq 2(|\om(C)| -1)< 2^\ell$,
 so we are done. Hence, assume that $k \le\ell-1$. By  Theorem~\ref{main},  $|\om(C)| \le \frac{2^k+1}{2^{k+1}}\cdot 2^\ell$, and thus
 $n \le \ell + k + \frac{2^k+1}{2^{k+1}}\cdot 2^\ell-1.$ If $\ell = 2$, then $k=1$, and so $p=\ell+k=3$ and $q\le\frac{2^k+1}{2^{k+1}}\cdot 2^\ell-1=2$, which is impossible.  Hence, $\ell\ge3$. Note that $k + \frac{2^k+1}{2^{k+1}}\cdot 2^\ell$ is maximized at $k=1$. Thus
 $k + \frac{2^k+1}{2^{k+1}}\cdot 2^\ell\le1+\frac34\cdot 2^\ell<2^\ell$ for
 $\ell\ge3$. Therefore, $n <2^\ell+ \ell-1$, which completes the proof.
}\end{proof}

\subsection{Cobipartite graphs}
 As the second application of Theorem~\ref{main}, we determine the maximum order of coreduced cobipartite graphs with a given corank and characterize the graphs achieving the maximum order.

From known relations between ranks of matrix sums (see the item 0.4.5\,(d) in \cite[p.~13]{hj}), we obtain the following:
\begin{lem}\label{lem:M+J}
For a symmetric matrix $M$,  $\rank(M+J)=\rank(M)+1$ if and only if $\1\notin \row(M)$.
\end{lem}
The following lemma is crucial for the proof of the main result of this section.
\begin{lem}\label{lem:subBB}

Let $B$ be a $p\times q$ $(0,1)$-matrix with $p\le q$, $\rank(B)=\ell$ and $\1\in\row(B)$.
Also assume that $B$ has no two identical columns or rows nor a zero row.
If $p+q\ge2^{\ell-1}+\ell-1$ and $\ell\ge6$, then $B$ is a submatrix of
\begin{equation}\label{eq:BB}
\left[
\begin{array}{c}
\BB_{\ell-1}\\ \hline  \1\\ \hline
J-\BB_{\ell-1}
\end{array}
\right],
\end{equation}
with a single exception in the case that $\ell=6$,  $p+q=2^{\ell-1}+\ell-1$, and the columns of $B$ are generated by
 \begin{equation}\label{eq:ell=6}
\left[
\begin{array}{c}
I_6\\ \hline  \x\\ \hline  \1\\ \hline
J_6-I_6\\ \hline
 \1-\x
\end{array}
\right],
\end{equation}
for some vector $\x$ of weight $2$ or $3$.
\end{lem}
\begin{proof}
We first construct a new matrix from $B$ as follows: if $\1$ is not already a row of $B$, we add it to the rows. Additionally, for any row $\x\ne\1$  of $B$, if $\1-\x$ is not a row, we add that as well. We call the resulting matrix $B'$.
 The matrix $B'$ is of the following form:
$$B'=\left[
\begin{array}{c}
B_0\\ \hline   \1\\ \hline
J-B_0
\end{array}
\right],$$
where $B_0$ consists of the rows of $B'$ whose first component is zero.
As $B'$ is obtained by adding some rows to $B$, it follows that $\rank(B')\ge\rank(B)$. However, each row of $B'$ can be expressed as a linear combination of $\1$ and some row of $B$. Since $\1 \in \row(B)$, this implies $\row(B')\subseteq\row(B)$, leading to $\rank(B')=\rank(B)=\ell$.
Given that $\1\notin\row(B_0)$ and every row of $B'$ can be formed through a linear combination of the rows of $B_0$ and $\1$, we conclude that $\rank(B_0)=\ell-1$. Our assumption on $B$ guarantees that $B_0$ has no two identical columns/rows and no zero rows.
If  $B_0$ has $\ell-1$ rows, then $B_0$ is a submatrix of $\BB_{\ell-1}$, and we are done.
Therefore, assume that $B_0$ has $\ell-1+k$ rows for some $k\ge1$.  So, $p\le2\ell+ 2k-1$.
By performing column-elementary operations and possibly permuting the rows, we can assume that $B_0$ has a basis of the form
$$\left[
\begin{array}{c}
I_{\ell-1}\\ \hline
C_{k\times (\ell-1)}
\end{array}
\right].$$
This basis has no identical rows nor a zero row. This implies that $C$ is a reduced matrix.
Every column of $B$ belongs to $\{A\bfb^\top:\bfb\in\om(C)\}$.  So $q\le|\om(C)|$.
 If $k\ge \ell-1$, then by Theorem~\ref{main}, $|\om(C)|\leq 2^{\ell-2}$. Thus
 $p+q\le2q \leq 2|\om(C)|\le 2^{\ell-1}$, which is a contradiction.
Hence, assume that $1\leq k\leq \ell-2 $.
By Theorem~\ref{main}, we have
$|\om(C)|\leq \frac{2^k+1}{2^{k+1}}\cdot 2^{\ell-1}$, and so
$$p+q \leq f:= 2\ell+ 2k-1+\frac{2^k+1}{2^{k+1}}\cdot 2^{\ell-1}.$$
If $\ell = 6$ and $2\leq k \leq 4$, by direct computation one can verify that $f < 2^{\ell-1}+\ell-1$.
For $\ell = 6$ and $k=1$, we have $f=2^{\ell-1}+\ell-1$. This implies that $q=|\om(C)|=\frac34\cdot 2^{\ell-1}$.
By the cases of equality  in Theorem~\ref{main}, $C$ should consists of a vector of weight $2$ or $3$, and thus $\col(B)$ has a basis of the form \eqref{eq:ell=6}.
If $\ell\ge7$,  $2k+ \frac{2^k+1}{2^{k+1}}\cdot 2^{\ell-1}$ is maximized at $k = 1$.
Therefore,
$$f\le 2\ell+1+\frac{3}{4}\cdot 2^{\ell-1} < 2^{\ell-1}+\ell-1,$$ from which the result follows.
\end{proof}

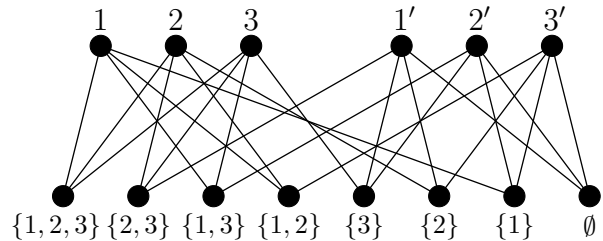
\begin{figure}[h!]\centering
\begin{tikzpicture}
\draw (-.46,.37)node {$3'$};
\draw (-1.46,.37)node {$2'$};
\draw (-2.46,.37)node {$1'$};
\draw (-4.5,.35)node {$3$};
\draw (-5.5,.35)node {$2$};
\draw (-6.5,.35)node {$1$};
\draw (0,-2.4)node {\footnotesize $\emptyset$};
\draw (-1,-2.4)node {\footnotesize $\{1\}$};
\draw (-2,-2.4)node {\footnotesize $\{2\}$};
\draw (-3,-2.4)node {\footnotesize $\{3\}$};
\draw (-4,-2.4)node {\footnotesize $\{1,2\}$};
\draw (-5,-2.4)node {\footnotesize $\{1,3\}$};
\draw (-6,-2.4)node {\footnotesize $\{2,3\}$};
\draw (-7.1,-2.4)node {\footnotesize $\{1,2,3\}$};
\draw [line width=.5pt] (-6.5,0)-- (-7,-2);
\draw [line width=.5pt] (-6.5,0)-- (-5,-2);
\draw [line width=.5pt] (-6.5,0)-- (-4,-2);
\draw [line width=.5pt] (-6.5,0)-- (-1,-2);
\draw [line width=.5pt] (-5.5,0)-- (-7,-2);
\draw [line width=.5pt] (-5.5,0)-- (-6,-2);
\draw [line width=.5pt] (-5.5,0)-- (-4,-2);
\draw [line width=.5pt] (-5.5,0)-- (-2,-2);
\draw [line width=.5pt] (-4.5,0)-- (-7,-2);
\draw [line width=.5pt] (-4.5,0)-- (-6,-2);
\draw [line width=.5pt] (-4.5,0)-- (-5,-2);
\draw [line width=.5pt] (-4.5,0)-- (-3,-2);
\draw [line width=.5pt] (-2.5,0)-- (-6,-2);
\draw [line width=.5pt] (-2.5,0)-- (-3,-2);
\draw [line width=.5pt] (-2.5,0)-- (-2,-2);
\draw [line width=.5pt] (-2.5,0)-- (0,-2);
\draw [line width=.5pt] (-1.5,0)-- (-5,-2);
\draw [line width=.5pt] (-1.5,0)-- (-3,-2);
\draw [line width=.5pt] (-1.5,0)-- (-1,-2);
\draw [line width=.5pt] (-1.5,0)-- (0,-2);
\draw [line width=.5pt] (-.5,0)-- (-4,-2);
\draw [line width=.5pt] (-.5,0)-- (-2,-2);
\draw [line width=.5pt] (-.5,0)-- (-1,-2);
\draw [line width=.5pt] (-.5,0)-- (0,-2);
\draw [fill=black] (-.5,0) circle (4pt);
\draw [fill=black] (-1.5,0) circle (4pt);
\draw [fill=black] (-2.5,0) circle (4pt);
\draw [fill=black] (-4.5,0) circle (4pt);
\draw [fill=black] (-5.5,0) circle (4pt);
\draw [fill=black] (-6.5,0) circle (4pt);
\draw [fill=black] (0,-2) circle (4pt);
\draw [fill=black] (-1,-2) circle (4pt);
\draw [fill=black] (-2,-2) circle (4pt);
\draw [fill=black] (-3,-2) circle (4pt);
\draw [fill=black] (-4,-2) circle (4pt);
\draw [fill=black] (-5,-2) circle (4pt);
\draw [fill=black] (-6,-2) circle (4pt);
\draw [fill=black] (-7,-2) circle (4pt);
\end{tikzpicture}
\caption{ The graph $\mathcal{D}_3$.}
\label{D6}
\end{figure}

We denote the bipartite graph $G$ with
$$B(G)=\left[
\begin{array}{c}
\BB_\ell\\ \hline
J-\BB_\ell
\end{array}
\right],$$
by $\mathcal{D}_\ell$. In other words, $\mathcal{D}_\ell$ is a bipartite graph with parts
 $\{1,1', \ldots, \ell,\ell'\}$ and ${\cal P}([\ell])$, such that each $S\in{\cal P}([\ell])$ has the $\ell$ neighbors $\{i: i\in  S\}\cup\{j':j\in [\ell]\setminus S\}$.  As an instance, $\mathcal{D}_3$ is depicted in Figure~\ref{D6}.

Now, we are in a position to prove the main result of this section. Recall that the complement of a graph $G$ is denoted by $\ov G$.
\begin{thm}
If $\ov G$ is a coreduced cobipartite graph with order $n$  and corank $r$, then
$$
n \leq
\begin{cases}
2^{\frac{r}{2}-1}+ r-2 &\ \mbox {$r$ even,}\\
2^{\frac{r-1}{2}}+ \frac{r-1}{2} &\ \mbox {$r$ odd.}
\end{cases}$$
The equality holds if and only if $G$ is isomorphic to ${\cal D}_{\frac{r}{2}-1}$ for even $r$, and to $\B_\frac{r-1}{2}$ for odd $r$.
\end{thm}
\begin{proof} Suppose that $\ov G$ is a coreduced cobipartite graph with corank $r$ and the maximum possible order $n$.
Let $\ov A=A(\ov G)$ and $A=A(G)$. Also let $B=B(G)$ be a $p\times q$  matrix. So, $n = p+q$. With no loss of generality, assume that $p\leq q$.
 Since $\ov G$ is a coreduced graph,  $G$ has no twins. So $B$ has no identical rows/columns.
  Note that $G$ might have an isolated vertex. In which case, we can assume that the isolated vertex lies in the larger part of $G$, that is, $B$ has a zero column rather than a zero row.
  Recall that $r=\rank(\ov A+I)$.
  So from $\ov A+I=J-A$, it follows that
\begin{equation}\label{eq:rankovA}
r-1\le\rank(A)=2\,\rank(B)\le r+1.
\end{equation}

We verified the result for $r\le10$ by a computer search. This is done by implementing an algorithm from \cite{ellin} (see also \cite{ack}) for constructing coreduced graphs of a fixed corank $r$.
For a given $r$, the input of the algorithm is the set of coreduced graphs with both order and corank equal to $r$ (which was generated by using McKay database of small graphs \cite{mk})
and the output of the algorithm is the set of all coreduced graphs of corank $r$.
So in what follows, we assume that $r\ge11$.

First suppose that $r=2\ell$ is even and so $\ell\ge6$. From \eqref{eq:rankovA} it follows that
 $\rank(A) = r$. Hence, by Lemma~\ref{lem:M+J}, $\1 \in \row(A)$.
It follows that $\1_q \in \row(B)$ and $\1^\top_p\in\col(B)$.
If $n=p+q<2^{\ell-1}+2\ell-2$, there is nothing to prove. Hence, we assume that $p+q\ge2^{\ell-1}+2\ell-2$.
 So $B$ satisfies the conditions of Lemma~\ref{lem:subBB}, and thus it is a submatrix of the matrix $C$ given in \eqref{eq:BB}.
 However, $\1^\top\notin\col(C)$ because $\col(C)$ has the following basis:
\begin{equation}\label{eq:ColBasis}
 \left[
\begin{array}{c|c}
\0^\top&I_{\ell-1}\\ \hline
1& \1_{\ell-1}\\ \hline
\1^\top&J_{\ell-1}-I_{\ell-1}
\end{array}
\right],
\end{equation}
and it is clear that such a basis cannot generate $\1^\top$.
Therefore, $B$ must have at least one row or one column less than $C$. This shows that $n\le2^{\ell-1}+2\ell-2$.
If we remove the $\1$ row of $C$, then the resulting matrix is $B({\cal D}_{\ell-1})$. So $G={\cal D}_{\ell-1}$, as desired.
To finish the proof, we show that if one deletes any other row or any column from $C$, then $\1^\top$ does not belong to the column space of the resulting matrix.
If we remove a row other than $\1$ from $C$ to obtain $C'$, then the restriction of \eqref{eq:ColBasis} to $C'$ forms a basis for $\col(C')$.
Again such a basis does not generate $\1^\top$.  A similar argument works in the case that $C'$ is obtained by removing one column from $C$.

Next, suppose that $r=2\ell-1$ is odd and so
$\ell\ge6$. Let $n \ge2^{\ell-1} + \ell-1$.
 To establish the theorem, it suffices to show that $G$ is isomorphic to $\B_{\ell-1}$.
 By \eqref{eq:rankovA}, we have
 $\rank(A)=2\ell-2$ or $2\ell$.
  If $\rank(A)=2\ell-2$, then
  we have necessarily $B=\BB_{\ell-1}$, that is $G=\B_{\ell-1}$ and
  we are done. So in what follows, we assume that
 $\rank(A)=2\ell$, i.e.
$\rank(B)=\ell$. 
Given that $A=J-(\ov A+I)$, we have $\rank(J-(\ov A+I))=r+1=\rank(-(\ov A+I))+1$. By invoking Lemma~\ref{lem:M+J}, this implies that $\1\notin\row(\ov A+I)$. Furthermore, since $\rank(J-A)<\rank(-A)$, another application of Lemma~\ref{lem:M+J} establishes that $\1\in\row(A)$, implying $\1 \in \row(B)$. Given this and the condition $n \geq 2^{\ell-1} + \ell-1$, the criteria outlined in Lemma~\ref{lem:subBB} are satisfied. Consequently, $\col(B)$ has a basis of the form \eqref{eq:ell=6} or  $B$ is a submatrix of \eqref{eq:BB}.
If the former occurs, then $\1^\top \notin \col(B)$, which implies $\1 \notin \row(A)$, leading to a contradiction. Therefore, $B$ is a submatrix of \eqref{eq:BB}.
Note that $\1$ cannot be a row of $B$. Otherwise, similar to the case of even $r$, we observe that $\1^\top \notin \col(B)$, resulting in $\1 \notin \row(A)$, which is a contradiction. 
Now, we make use of the fact that $\1\notin\row(\ov A+I)$. We have
$$\ov A+I=
\left[
\begin{array}{c|c}
J&J-B\\ \hline
J-B^\top&J
\end{array}
\right].$$
We claim that if some vector $\x$ is a row of $B$, then $\1-\x$ is not a row of $B$. If this fails, then we can obtain
$\left[\begin{array}{c|c}2\1_p&\1_q\end{array}\right]$ as sum of two rows of $\left[\begin{array}{c|c}J&J-B\end{array}\right]$.
Also, as $B$ has more than
$2^{\ell-2}$ columns, it contains some two columns of the forms $\y^\top$ and $\1^\top-\y^\top$.
The two corresponding rows  in
$\left[\begin{array}{c|c}J-B^\top&J\end{array}\right]$
 sum up to
$\left[\begin{array}{c|c}\1_p&2\1_q\end{array}\right]$. It turns out that $\1_n=\frac13\left[\begin{array}{c|c}2\1_p&\1_q\end{array}\right]+\frac13\left[\begin{array}{c|c}\1_p&2\1_q\end{array}\right]\in\row(\ov A+I)$, again a contradiction.
This proves the claim.
So we have established that $B$ is a submatrix of \eqref{eq:BB} such that
 $\1_q$  is not a row of $B$ and if $\x$ is a row of $B$, then  $\1-\x$ is not a row of $B$. It follows that $B$ has at most $\ell-1$ rows.
This is a contradiction because $\rank(B)=\ell$. This means that
the case $\rank(A)=2\ell$ is impossible, and the proof is complete.
 \end{proof}

\section*{Acknowledgments}
The second author carried out part of this work during a Humboldt Research Fellowship at the University of Hamburg. He thanks the Alexander von Humboldt-Stiftung for financial support. The authors thank Omid Etesami for fruitful discussions, particularly for the proof of Lemma~\ref{lem:Omega(v)}. They also thank an anonymous referee for several helpful comments.

{}
\end{document}